\documentclass[12pt]{article}

\usepackage{graphicx}

\def\be{\begin{equation}}
\def\ee{\end{equation}}
\newcommand{\kk}[2]{\frac{#1}{#2}}

\newcommand{\ff}[1]{{\bf  #1}}
\def\a{\alpha}
\def\b{\beta}

\def\e{\epsilon}

\def\s{\qquad}
\def\ra{\rightarrow}
\def\={\approx}
\def\8{{\infty}}
\def\x{\ff{x}}
\def\vcode#1#2#3#4{\begin{figure}
\begin{center}
\begin{minipage}[c]{#1\textwidth}
{{\small #2 \hrule \vspace{5pt}   %
{\it #3}  \vspace{5pt} \hrule }}
\end{minipage}
\caption{#4}
\end{center}   \end{figure}} %%
\begin{document}

\title{Firefly Algorithms for Multimodal Optimization}

\author{Xin-She Yang \\
Department of Engineering, University of Cambridge, \\
Trumpington Street, Cambridge CB2 1PZ, UK }

\date{}

\maketitle

%% Abstract

\abstract{
Nature-inspired algorithms are among the most powerful algorithms for
optimization.  This paper intends to provide a detailed description of
a new Firefly Algorithm (FA) for multimodal optimization applications. We will
compare the proposed firefly algorithm with other metaheuristic
algorithms such as particle swarm optimization (PSO).
Simulations and results indicate that the proposed firefly algorithm is
superior to existing metaheuristic algorithms.
Finally we will discuss its  applications and implications
for further research. \\
}

\noindent {\bf Citation detail:}
X.-S. Yang, ``Firefly algorithms for multimodal optimization", in: {\it Stochastic Algorithms:
Foundations and Applications}, SAGA 2009, Lecture Notes in Computer Sciences, Vol. {\bf 5792}, pp. 169-178 (2009).

%% Begin of Main Text %%

\section{Introduction}

Biologically inspired algorithms are becoming powerful in modern
numerical optimization \cite{Bod,Deb,Gold,Ken2,Yang,Yang2},
especially for the NP-hard problems such
as the travelling salesman problem. Among these biology-derived algorithms,
the multi-agent metaheuristic algorithms such as particle swarm optimization
form hot research topics in the start-of-the-art algorithm development in
optimization and other applications \cite{Bod,Deb,Yang}.

Particle swarm optimization (PSO) was developed by Kennedy and
Eberhart in 1995 \cite{Ken}, based on the swarm behaviour such
as fish and bird schooling in nature, the so-called swarm intelligence.
Though particle swarm optimization has many
similarities with genetic algorithms, but
it is much simpler because it does not use mutation/crossover
operators. Instead, it uses the real-number randomness
and the global communication among the swarming particles. In this
sense, it is also easier to implement as it uses mainly real numbers.

This paper aims to introduce the new Firefly Algorithm and to provide
the comparison study of the FA with PSO and other relevant algorithms.
We will first outline the particle swarm optimization, then formulate
the firefly algorithms and finally give the comparison about the performance
of these algorithms. The FA optimization
seems more promising than particle swarm optimization in the sense
that FA can deal with multimodal functions more naturally and efficiently.
In addition, particle swarm optimization is just a special class of
the firefly algorithms as we will demonstrate this in this paper.

\section{Particle Swarm Optimization}

\subsection{Standard PSO}

The PSO algorithm searches the space of the objective functions
by adjusting the trajectories of individual agents,
called particles, as the piecewise paths formed by positional
vectors in a quasi-stochastic manner \cite{Ken,Ken2}.
There are now as many as about 20 different variants of PSO.
Here we only describe the simplest and yet popular standard PSO.

The particle movement has two
major components: a stochastic component and a deterministic component.
A particle is attracted toward the position of the current global best
$\ff{g}^*$ and its own best location $\x_i^*$ in history,
while at the same time it has a tendency to move randomly.
When a particle finds a location that is better than any previously found locations,
then it updates it as the new current best for particle $i$. There is
a current global best for all $n$ particles. The aim is to find the global best
among all the current best solutions until the objective no longer improves or
after a certain number of iterations.

For the particle movement, we use $\x^*_i$ to denote the current best for
particle $i$, and $\ff{g}^* \=\min$ or $\max \{f(\x_i)\} (i=1,2,...,n)$ to denote
the current global best. Let $\x_i$ and $\ff{v}_i$ be the position vector and velocity for
particle $i$, respectively. The new velocity vector is determined by the
following formula
\be \ff{v}_i^{t+1}= \ff{v}_i^t  + \a \ff{\e}_1  \odot
(\ff{g}^*-\x_i^t) + \b \ff{\e}_2 \odot (\x_i^*-\x_i^t).
\label{pso-speed-100}
\ee
where $\ff{\e}_1$ and $\ff{\e}_2$ are two random vectors, and each
entry taking the values between 0 and 1. The Hadamard product of two
matrices $\ff{u \odot v}$ is defined as the entrywise product, that
is $[\ff{u \odot v}]_{ij}=u_{ij} v_{ij}$.
The parameters $\a$ and $\b$ are the learning parameters or
acceleration constants, which can typically be taken as, say, $\a\=\b \=2$.
The initial values of $\x_i^{t=0}$ can be taken as the bounds or limits $a=\min (x_j)$,
$b=\max(x_j)$ and $\ff{v}_i^{t=0}\!\!=0$.
The new position can then be updated by \be
\x_i^{t+1}=\x_i^t+\ff{v}_i^{t+1}. \ee Although $\ff{v}_i$ can be any
values, it is usually bounded in some range $[0, \ff{v}_{max}]$.

There are many variants which extend the standard PSO
algorithm,
and the most noticeable improvement is probably to use inertia function $\theta
(t)$ so that $\ff{v}_i^t$ is replaced by $\theta(t) \ff{v}_i^t$
where $\theta$ takes the values between 0 and 1. In the simplest case,
the inertia function can be taken as a constant, typically $\theta \= 0.5 \sim 0.9$.
This is equivalent to introducing a virtual mass to stabilize the motion
of the particles, and thus the algorithm is expected to converge more quickly.

\section{Firefly Algorithm}

\subsection{Behaviour of Fireflies}

The flashing light of fireflies is an amazing sight in the summer sky in
the tropical and temperate regions. There are about two thousand
firefly species, and most fireflies produce short and rhythmic flashes.
The pattern of flashes is often unique for a particular species.
The flashing light is produced by a process of bioluminescence, and
the true functions of such signaling systems are still debating. However,
two fundamental functions of such flashes are to attract mating partners
(communication), and to attract potential prey. In addition,
flashing may also serve as a protective warning mechanism.
The rhythmic flash, the rate of flashing and the amount of time
form part of the signal system that brings both sexes together.
Females respond to a male's unique pattern of flashing
in the same species, while in some species such as photuris, female
fireflies  can mimic the mating
flashing pattern of other species so as to lure and eat the male
fireflies who may mistake the flashes as a potential suitable mate.

We know that the light intensity at a particular distance $r$ from
the light source obeys the inverse square law. That is to say, the light
intensity $I$ decreases as the distance $r$ increases in terms of
$I \propto 1/r^2$. Furthermore, the air absorbs light which becomes
weaker and weaker as the distance increases. These two combined factors
make most fireflies visible only to a limited distance, usually several hundred
meters at night, which is usually good enough for fireflies to communicate.

The flashing light can be formulated in such a way that it is associated with
the objective function to be optimized, which makes it possible to formulate
new optimization algorithms. In the rest of this paper, we will first outline
the basic formulation of the Firefly Algorithm (FA) and then discuss
the implementation as well as its analysis in detail.

\subsection{Firefly Algorithm}

Now we can idealize some of the flashing characteristics of fireflies
so as to develop firefly-inspired algorithms. For simplicity in
describing our new Fireflire Algorithm (FA), we now use the following
three idealized rules: 1) all fireflies are unisex so that one firefly will be attracted to other fireflies
regardless of their sex; 2) Attractiveness is proportional to their brightness,
thus for any two flashing fireflies, the less brighter one will move towards the brighter
one. The attractiveness is proportional to the brightness and they both
decrease as their distance increases. If there is
no brighter one than a particular firefly, it will move randomly;
3) The brightness of a firefly is affected or  determined by the landscape of the
objective function. For a  maximization problem, the brightness can simply be proportional
to the value of the objective function. Other forms of brightness can be defined in a similar
way to the fitness function in genetic algorithms.

Based on these three rules, the basic steps of the firefly algorithm (FA)
can be summarized as the pseudo code shown in Fig. \ref{fa-fig-100}.

\vcode{0.9}{{\sf Firefly Algorithm}} {
\indent \quad Objective function $f(\x), \s \x=(x_1, ..., x_d)^T$ \\
\indent \quad Generate initial population of fireflies $\x_i \; (i=1,2,...,n)$ \\
\indent \quad Light intensity $I_i$ at $\x_i$ is determined by $f(\x_i)$  \\
\indent \quad Define light absorption coefficient $\gamma$ \\
\indent \quad {\bf while} ($t<$MaxGeneration) \\
\indent \quad {\bf for} $i=1:n$ all $n$ fireflies \\
\indent \qquad {\bf for} $j=1:i$ all $n$ fireflies \\
\indent \qquad \qquad {\bf if} ($I_j>I_i$), Move firefly $i$ towards $j$ in d-dimension; {\bf end if} \\
\indent \qquad \qquad Attractiveness varies with distance $r$ via $\exp[-\gamma r]$ \\
\indent \qquad \qquad Evaluate new solutions and update light intensity \\
\indent \qquad {\bf end for }$j$ \\
\indent \quad {\bf end for }$i$ \\
\indent \quad Rank the fireflies and find the current best  \\
\indent \quad {\bf end while} \\
\indent \quad Postprocess results and visualization }{Pseudo code of the firefly algorithm (FA). \label{fa-fig-100} }

In certain sense, there is some conceptual similarity between the firefly algorithms and the bacterial
foraging algorithm (BFA) \cite{Gazi,Passino}. In BFA, the attraction among bacteria
is based partly on their fitness and partly on their distance, while in FA, the attractiveness is
linked to their objective function and monotonic decay of the attractiveness with distance.
However, the agents in FA have adjustable visibility and more versatile in attractiveness variations,
which usually leads to higher mobility and thus the search space is explored more efficiently.

\subsection{Attractiveness}

In the firefly algorithm, there are two important issues: the variation of
light intensity and formulation of the attractiveness.
For simplicity, we can always assume that the attractiveness
of a firefly is determined by its brightness which in turn is associated with
the encoded objective function.

In the simplest case for maximum optimization problems,
the brightness $I$ of a firefly at a particular location $\x$ can be chosen
as $I(\x) \propto f(\x)$. However, the attractiveness $\b$ is relative, it should be
seen in the eyes of the beholder or judged by the other fireflies. Thus, it will
vary with the distance $r_{ij}$ between firefly $i$ and firefly $j$. In addition,
light intensity decreases with the distance from its source, and light is also
absorbed in the media,
so we should allow the attractiveness to vary with the degree of absorption.
In the simplest form, the light intensity $I(r)$ varies according to the inverse
square law $I(r)=I_s/r^2$ where $I_s$ is the intensity at the source.
For a given medium with a fixed light absorption coefficient $\gamma$,
the light intensity $I$ varies with the distance $r$. That is
$ I=I_0 e^{-\gamma r}, $ where $I_0$ is the original light intensity.
In order to avoid the singularity
at $r=0$ in the expression $I_s/r^2$, the combined effect of both the inverse
square law and absorption can be approximated using the following Gaussian form
\be I(r)=I_0 e^{-\gamma r^2}. \ee
Sometimes, we may need a function which decreases monotonically at a slower rate. In
this case, we can use the following  approximation
\be I(r)=\kk{I_0}{1+\gamma r^2}.  \ee
At a shorter distance, the above two forms are essentially the same.
This is because the series expansions about $r=0$
\be e^{-\gamma r^2} \= 1- \gamma r^2 + \kk{1}{2} \gamma^2 r^4 + ..., \s
 \kk{1}{1+\gamma r^2} \= 1-\gamma r^2 + \gamma^2 r^4 + ..., \ee
are equivalent to each other up to the order of $O(r^3)$.

As a firefly's attractiveness is proportional to
the light intensity seen by adjacent fireflies, we can now define
the attractiveness $\b$ of a firefly by
\be \b (r)  = \b_0 e^{-\gamma r^2}, \label{att-equ-100} \ee
where $\b_0$ is the attractiveness at $r=0$. As it is often faster to calculate
$1/(1+r^2)$ than an exponential function, the above function, if necessary, can
conveniently be replaced by $ \b = \kk{\b_0}{1+\gamma r^2}$.
Equation (\ref{att-equ-100}) defines a characteristic distance
$\Gamma=1/\sqrt{\gamma}$ over which the attractiveness changes significantly
from $\b_0$ to $\b_0 e^{-1}$.

In the implementation, the actual form of attractiveness function $\b(r)$ can be
any monotonically decreasing functions such as the following generalized form
\be \b(r) =\b_0 e^{-\gamma r^m}, \s (m \ge 1). \ee
For a fixed $\gamma$, the characteristic length becomes
$\Gamma=\gamma^{-1/m} \ra 1$ as $m \ra \infty$.  Conversely,
for a given length scale $\Gamma$ in an optimization problem, the parameter $\gamma$
can be used as a typical initial value. That is $\gamma =\kk{1}{\Gamma^m}$.

\subsection{Distance and Movement}

The distance between any two fireflies $i$ and $j$ at $\x_i$ and $\x_j$, respectively, is
the Cartesian distance
\be r_{ij}=||\x_i-\x_j|| =\sqrt{\sum_{k=1}^d (x_{i,k} - x_{j,k})^2}, \ee
where $x_{i,k}$ is the $k$th component of the spatial coordinate $\x_i$ of $i$th
firefly. In 2-D case, we have $r_{ij}=\sqrt{(x_i-x_j)^2+(y_i-y_j)^2}.$

The movement of a firefly $i$ is attracted to another more attractive (brighter)
firefly $j$ is determined by
\be \x_i =\x_i + \b_0 e^{-\gamma r^2_{ij}} (\x_j-\x_i) + \a \; ({\rm rand}-\kk{1}{2}), \ee
where the second term is due to the attraction while the third term
is randomization with $\a$ being the randomization parameter. ${\sf rand}$ is a
random number generator uniformly distributed in $[0,1]$. For most cases in our implementation,
we can take $\b_0=1$ and $\a \in [0,1]$. Furthermore, the randomization term can easily be
extended to a normal distribution $N(0,1)$ or other distributions.
In addition, if the scales vary significantly in different dimensions such as $-10^5$ to $10^5$ in
one dimension while, say, $-0.001$ to $0.01$ along the other, it is a good idea to
replace $\a$ by $\a S_k$ where the scaling parameters $S_k (k=1,...,d)$ in
the $d$ dimensions should be determined by the actual scales of the problem of interest.

The parameter $\gamma$ now characterizes the variation of the attractiveness,
and its value is crucially important in determining the speed of the convergence
and how the FA algorithm behaves. In theory, $\gamma \in [0,\infty)$, but in
practice, $\gamma=O(1)$ is determined by the characteristic length $\Gamma$ of the
system to be optimized. Thus, in most applications, it typically varies
from $0.01$ to $100$.

\subsection{Scaling and Asymptotic Cases}

It is worth pointing out that the distance $r$ defined above is {\it not} limited
to the Euclidean distance. We can define many other forms of distance $r$ in the $n$-dimensional
hyperspace, depending on the type of problem of our interest. For example,
for job scheduling problems, $r$ can be defined as the time lag or time interval.
For complicated networks such as the Internet and social networks, the distance
$r$ can be defined as the combination of the degree of local clustering and
the average proximity of vertices. In fact, any measure that can effectively characterize
the quantities of interest in the optimization problem can be used as the `distance' $r$.
The typical scale $\Gamma$ should be associated with the scale in
the optimization problem of interest. If $\Gamma$ is the typical scale for a given
optimization problem, for a very large number of fireflies $n \gg m$ where
$m$ is the number of local optima, then the initial locations of these $n$
fireflies should distribute relatively uniformly over
the entire search space in a similar manner as the initialization of
quasi-Monte Carlo simulations.
As the iterations proceed, the fireflies would converge
into all the local optima (including the global ones) in a stochastic manner.
By comparing the best
solutions among all these optima, the global optima can easily be achieved.
At the moment, we are trying to formally prove that the firefly algorithm
will approach global optima when $n \ra \infty$ and $t \gg 1$. In reality,
it converges very quickly, typically with less than 50 to 100 generations,
and this will be demonstrated using various standard test functions later in this paper.

There are two important limiting cases when $\gamma \ra 0$ and $\gamma \ra \infty$.
For $\gamma \ra 0$, the attractiveness is constant $\b=\b_0$ and $\Gamma \ra \infty$,
this is equivalent to say that the light intensity does not decrease in an idealized sky.
Thus, a flashing firefly can be seen anywhere in the domain. Thus, a single (usually global)
optimum can easily be reached. This corresponds to a special case of particle
swarm optimization (PSO) discussed earlier. Subsequently, the efficiency of this special case
is the same as that of PSO.

On the other hand, the limiting case $\gamma \ra \infty$ leads to
$\Gamma \ra 0$ and  $\b(r) \ra \delta(r)$ (the Dirac delta function), which means that
the attractiveness is almost zero in the sight of other fireflies or the fireflies are short-sighted. This is equivalent
to the case where the fireflies fly in a very foggy region randomly. No other fireflies can be
seen, and each firefly roams in a completely random way. Therefore, this corresponds
to the completely random search method.
As the firefly algorithm is usually in somewhere between these two extremes, it is possible
to adjust the parameter $\gamma$ and $\alpha$ so that it can outperform both the random search
and PSO. In fact, FA can find the global optima as well as all the local optima
simultaneously in a very effective manner.
This advantage will be demonstrated in detail later in the implementation.
A further advantage of FA is that different fireflies will work almost independently, it is
thus particularly suitable for parallel implementation. It is even better than genetic algorithms
and PSO because fireflies aggregate more closely around each optimum (without jumping around
as in the case of genetic algorithms). The interactions between different subregions are
minimal in parallel implementation.
\begin{figure}
\vspace{-15pt}
\centerline{\includegraphics[width=2.5in,height=2in]{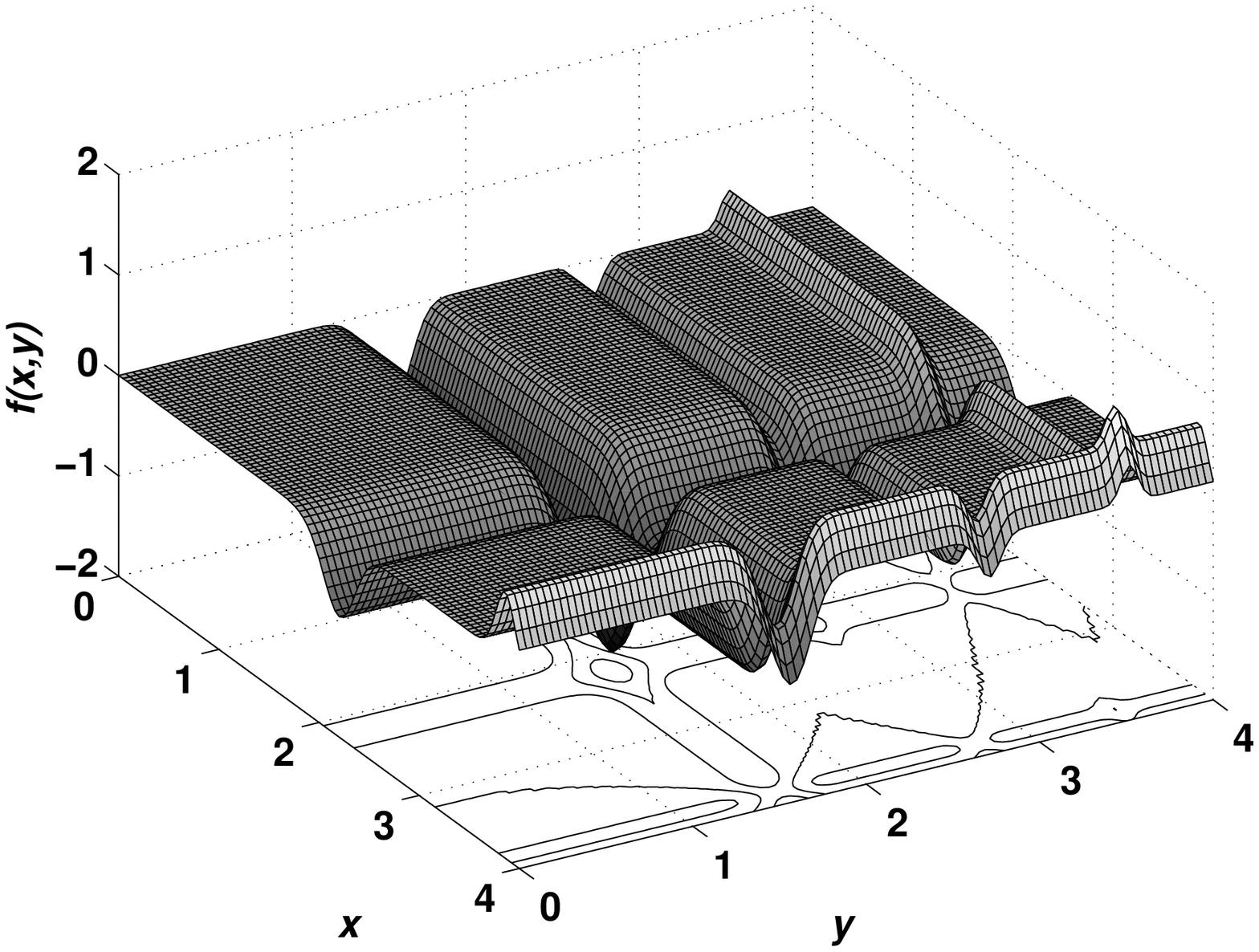} }
\vspace{-10pt}
\caption{Michalewicz's function for two independent variables
with a global minimum $f_* \approx -1.801$ at $(2.20319,1.57049)$. \label{yangfig-100} }
\end{figure}

\section{Multimodal Optimization with Multiple Optima}

\subsection{Validation}

\begin{figure}
 \centerline{\includegraphics[height=1.7in,width=2in]{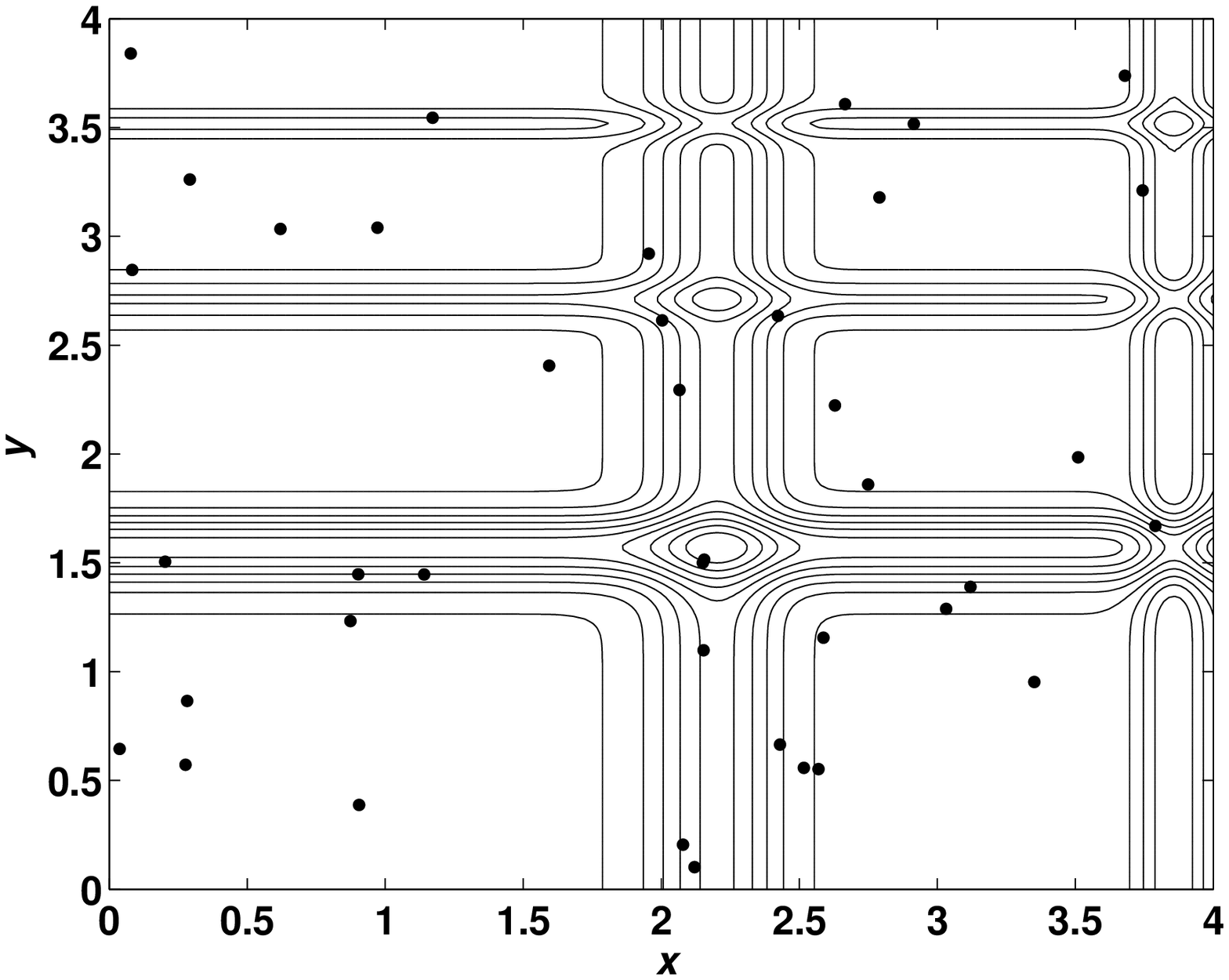}
 \includegraphics[height=1.7in,width=2in]{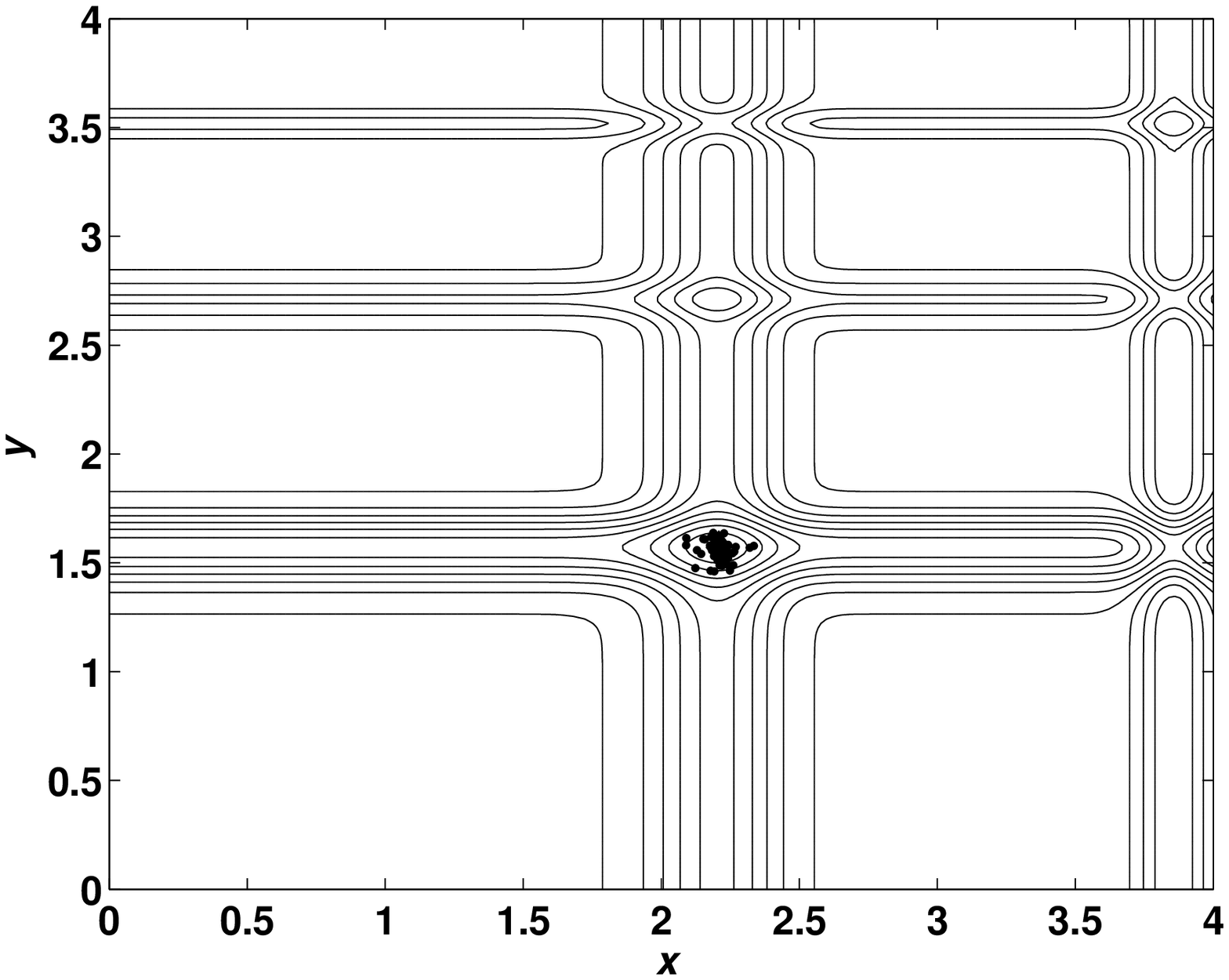}}
\vspace{-5pt}
\caption{The initial 40 fireflies (left) and their
locations after $10$ iterations (right). \label{yangfig-200} }
\end{figure}

In order to demonstrate how the firefly algorithm works, we have implemented it in
Matlab. We will use various test functions to validate the new algorithm.
As an example, we now use the FA to find the
global optimum of the Michalewicz function
\be f(\x)=-\sum_{i=1}^d \sin (x_i) [\sin (\kk{i x_i^2}{\pi})]^{2m}, \ee
where $m=10$ and $d=1,2,...$. The global minimum $f_* \approx -1.801$ in 2-D occurs
at $(2.20319,1.57049)$, which can be found after about 400 evaluations
for 40 fireflies after 10 iterations (see Fig. \ref{yangfig-100} and Fig. \ref{yangfig-200}).
Now let us use the FA to find the optima of some tougher
test functions. This is much more efficient than most of existing metaheuristic algorithms.
In the above simulations,  the values
of the parameters are $\a=0.2$, $\gamma=1$ and $\b_0=1$.

\begin{figure}
\vspace{-25pt}
\centerline{\includegraphics[width=3in,height=2.5in]{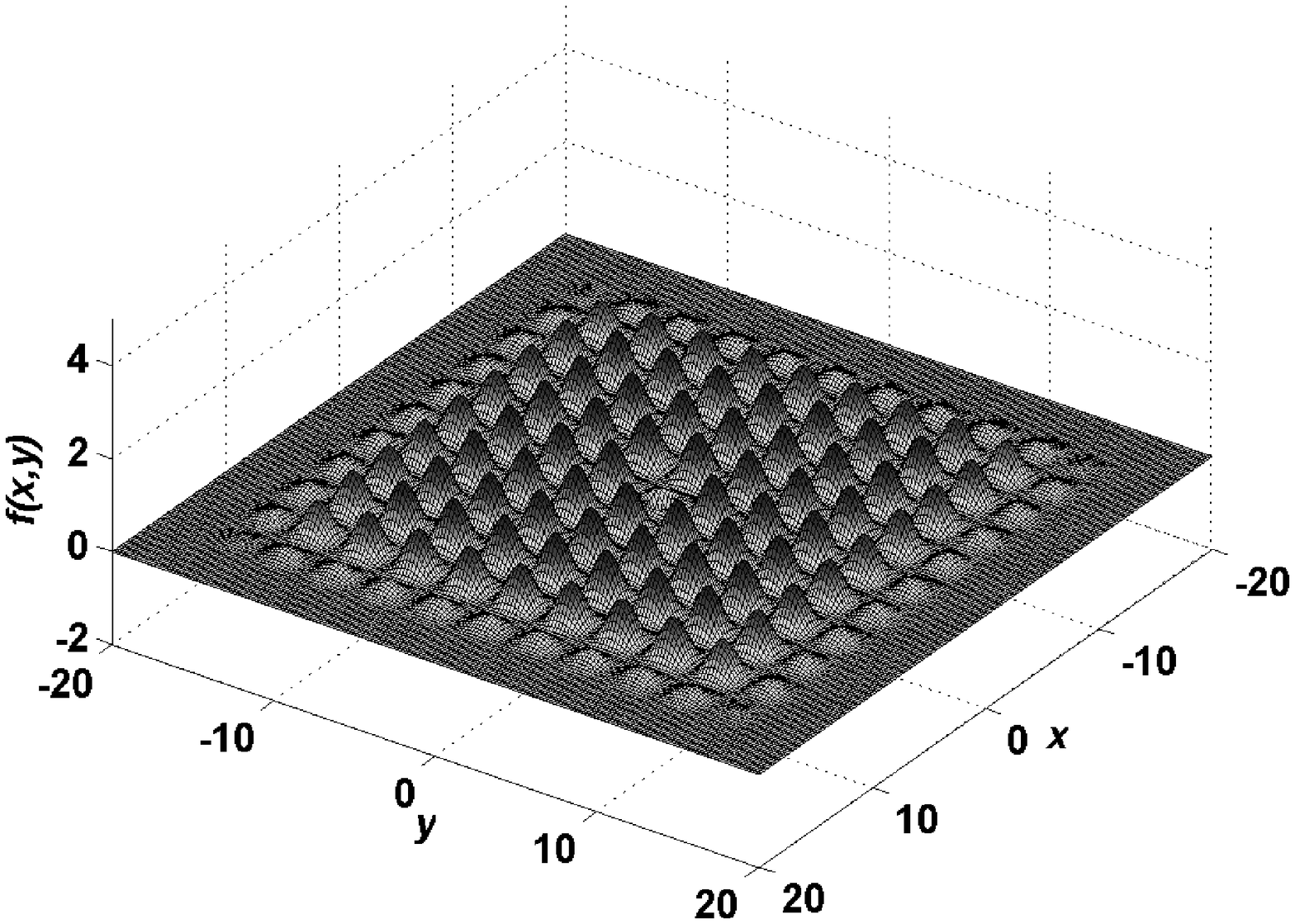} }
\vspace{-10pt}
\caption{Yang's function in 2D with a global minimum $f_* = -1$ at $(0,0)$ where $a=15$. \label{yangfig-300} }
\end{figure}

We have also used much tougher test functions. For example, Yang described a multimodal
function which looks like a standing-wave pattern  \cite{Yang3}
 \be f(\x)=\Big[ e^{-\sum_{i=1}^d (x_i/a)^{2 m}} - 2 e^{-\sum_{i=1}^d x_i^2} \Big] \cdot
 \prod_{i=1}^d \cos^2 x_i, \s m=5, \ee
is multimodal with many local peaks and valleys, and it has a unique
global minimum $f_*=-1$ at $(0,0,...,0)$ in the region $-20 \le x_i \le 20$ where
$i=1,2,...,d$ and $a=15$. The 2D landscape of Yang's function is shown in Fig. \ref{yangfig-300}.

\subsection{Comparison of FA with PSO and GA}

Various studies show that PSO algorithms can
outperform genetic algorithms (GA) \cite{Gold} and other conventional algorithms for
solving many optimization problems. This is
partially due to that fact that the broadcasting ability of the current
best estimates gives better and quicker convergence towards the
optimality. A general framework for evaluating statistical performance of evolutionary algorithms
has been discussed in detail by Shilane et al. \cite{Shilane}.
Now we will compare the Firefly Algorithms with PSO, and genetic
algorithms for various standard test functions.
For genetic algorithms, we have used the standard version with no elitism with
a mutation probability of $p_m=0.05$ and a crossover probability of $0.95$.
For the particle swarm optimization, we have also used the standard version with
the learning parameters $\a \= \b \= 2$ without the inertia correction \cite{Gold,Ken,Ken2}.
We have used various population sizes from $n=15$ to $200$, and found that for most
problems, it is sufficient to use $n=15$ to $50$. Therefore,
we have used a fixed population size of $n=40$ in all our simulations for comparison.

After implementing these algorithms using
Matlab, we have carried out extensive simulations and each algorithm has been
run at least 100 times so as to carry out meaningful statistical analysis.
The algorithms stop when the variations of function values are less than
a given tolerance $\e \le 10^{-5}$. The results are
summarized in the following table (see Table 1) where the global optima
are reached. The numbers are in the
format: average number of evaluations (success rate), so
$3752 \pm 725 (99\%)$ means that the average number (mean) of function
evaluations is 3752 with a standard deviation of 725. The success rate
of finding the global optima for this algorithm is $99\%$.

\begin{table}[ht]

\caption{Comparison of algorithm performance}
\centering

\begin{tabular}{ccccc}
\hline \hline
Functions/Algorithms & GA & PSO & FA \\
\hline

 Michalewicz's ($d\!\!=\!\!16$)  & $89325 \pm 7914 (95 \%)$  & $6922 \pm 537 (98\%)$  & $3752 \pm 725 (99\%)$ \\

 Rosenbrock's ($d\!\!=\!\!16$) & $55723 \pm 8901 (90\%)$ & $32756 \pm 5325 (98\%)$ & $7792 \pm 2923 (99\%) $ \\

 De Jong's ($d\!\!=\!\!256$) & $25412 \pm 1237 (100\%)$ & $17040 \pm 1123 (100\%)$ & $7217 \pm 730 (100\%)$\\
 Schwefel's ($d\!\!=\!\!128$) & $227329 \pm 7572 (95\%)$ & $14522 \pm 1275 (97\%)$ & $9902 \pm 592 (100\%)$ \\

 Ackley's ($d\!\!=\!\!128$) & $32720 \pm 3327 (90\%)$ & $23407 \pm 4325 (92\%)$ & $5293 \pm 4920 (100\%)$ \\

 Rastrigin's & $110523 \pm 5199 (77 \%)$ & $79491 \pm 3715 (90\%)$ & $15573 \pm 4399 (100\%)$ \\

 Easom's & $19239 \pm 3307 (92\%)$ & $17273 \pm 2929 (90\%)$ & $7925 \pm 1799 (100\%)$ \\

 Griewank's & $70925 \pm 7652 (90\%)$ & $55970 \pm 4223 (92\%)$ & $12592 \pm 3715 (100\%)$ \\

 Shubert's (18 minima) & $54077 \pm 4997 (89\%)$ & $23992 \pm 3755 (92\%)$ & $12577 \pm 2356 (100\%)$ \\

 Yang's ($d=16$) & $27923 \pm 3025 (83\%)$ & $14116 \pm 2949 (90\%)$ & $7390 \pm 2189 (100\%)$ \\

\hline
\end{tabular}
\end{table}

We can see that the FA is much more efficient in finding the global optima
with higher success rates. Each function evaluation is virtually instantaneous
on modern personal computer. For example, the computing time for 10,000 evaluations
on a 3GHz desktop is about 5 seconds. Even with graphics for displaying the
locations of the particles and fireflies, it usually takes less than a few minutes.
It is worth pointing out that more formal statistical hypothesis testing
can be used to verify such significance.

\section{Conclusions}

In this paper, we have formulated a new firefly algorithm and analyzed its
similarities and differences with particle swarm optimization. We then implemented
and compared these algorithms. Our simulation results for finding the global optima
of various test functions suggest that particle swarm often outperforms
traditional algorithms such as genetic algorithms, while the new firefly algorithm
is superior to both PSO and GA in terms of both efficiency and success rate.
This implies that FA is potentially more powerful in solving NP-hard problems
which will be investigated further in future studies.

The basic firefly algorithm is very efficient, but we can see that the solutions
are still changing as the optima are approaching. It is possible to improve the
solution quality by reducing the randomness gradually. A further
improvement on the convergence of the algorithm is to vary
the randomization parameter $\a$ so that it decreases gradually as the optima
are approaching. These could form important topics for further research.
Furthermore, as a relatively straightforward extension, the Firefly Algorithm can be modified
to solve multiobjective optimization problems.
In addition, the application of firefly algorithms in combination
with other algorithms may form an exciting area for further research.

%% End of text %%


\begin{thebibliography}{A}
\bibitem{Bod}  Bonabeau E., Dorigo M., Theraulaz G., {\it Swarm Intelligence:
From Natural to Artificial Systems}. Oxford University Press, (1999)


\bibitem{Deb} Deb. K., {\it Optimisation for Engineering Design}, Prentice-Hall, New Delhi, (1995).


\bibitem{Gazi} Gazi K., and Passino K. M., Stability analysis of social foraging swarms,
{\it IEEE Trans. Sys. Man. Cyber. Part B - Cybernetics}, {\bf 34},
539-557 (2004).

\bibitem{Gold} Goldberg D. E., {\it Genetic Algorithms in Search, Optimisation and
Machine Learning}, Reading, Mass.: Addison Wesley (1989).

\bibitem{Ken} Kennedy J. and Eberhart R. C.: Particle swarm optimization. {\it
Proc. of IEEE International Conference on Neural Networks},
Piscataway, NJ. pp. 1942-1948 (1995).

\bibitem{Ken2} Kennedy J., Eberhart R., Shi Y.: {\it Swarm intelligence}, Academic
Press, (2001).

\bibitem{Passino} Passino K. M., {\it Biomimicrt of Bacterial Foraging for Distributed
Optimization}, University Press, Princeton, New Jersey (2001).


\bibitem{Shilane} Shilane D., Martikainen J., Dudoit S., Ovaska S. J.,
 A general framework for statistical performance comparison
of evolutionary computation algorithms, {\it Information Sciences: an Int. Journal},
{\bf 178}, 2870-2879 (2008).

\bibitem{Yang} Yang X. S., {\it Nature-Inspired Metaheuristic Algorithms},
Luniver Press, (2008).

\bibitem{Yang2} Yang X. S., Biology-derived algorithms in engineering optimizaton
(Chapter 32), in {\it Handbook of Bioinspired Algorithms and Applications}
(eds Olarius \& Zomaya),  Chapman \& Hall / CRC (2005).

\bibitem{Yang3} Yang X. S., {\it Engineering Optimization: An Introduction with Metaheuristic Applications},
Wiley (2010).


\end{thebibliography}
\end{document}